\documentclass[a4paper]{article}

\usepackage[english]{babel}
\usepackage[utf8x]{inputenc}
\usepackage[T1]{fontenc}
\usepackage{lipsum}
\usepackage{blindtext}
\usepackage{graphicx,amssymb,amsmath,textcomp}
\usepackage{newpxtext} 

\usepackage[a4paper,top=3cm,bottom=2cm,left=3cm,right=3cm,marginparwidth=1.75cm]{geometry}
\usepackage{tikz}
\usetikzlibrary{matrix,arrows,positioning,calc}
\usepackage{verbatim}

\usepackage{amsmath,amsthm}
\usepackage{faktor}
\usepackage{xfrac} 
\usepackage{graphicx}
\usepackage[colorinlistoftodos]{todonotes}
\usepackage[colorlinks=true, allcolors=blue]{hyperref}
\usepackage{mathtools}

\title{\textbf{Notes on Derived Geometric Formulations in Physics }}
\author{Kadri İlker Berktav \footnote{University of Zurich, Institute of Mathematics, Zurich, Switzerland; e-mail: kadriilker.berktav@math.uzh.ch } \\ 
} 

\date{\vspace{-5ex}}

\begin{document}
\theoremstyle{plain}
\newtheorem{theorem}{Theorem}[section] 
\newtheorem{lemma}[theorem]{Lemma} 
\newtheorem{proposition}[theorem]{Proposition} 
\newtheorem{corollary}[theorem]{Corollary}
\newtheorem{conjecture}[theorem]{Conjecture}
\theoremstyle{definition}
\newtheorem{notations}[theorem]{Notations}
\newtheorem{notation}[theorem]{Notation}
\newtheorem{remark}[theorem]{Remark}
\newtheorem{observation}[theorem]{Observation}
\newtheorem{definition}[theorem]{Definition}
\newtheorem{condition}[theorem]{Condition}
\newtheorem{example}[theorem]{Example}
\let\pf\proof
\let\epf\endproof
\numberwithin{equation}{section} 

\maketitle

\begin{abstract}
	This is an overview of  higher structural constructions in physics. The main motivations of our current attempt are as follows: \textit{(i)} to provide a brief introduction to  \textit{derived algebraic geometry}, \textit{(ii)} to understand how  \textit{derived objects} naturally appear in physics and give rise to  formal mathematical treatment, and \textit{(iii)} to investigate how  \textit{factorization algebras} together with certain higher categorical structures come into play to encode the structure of \textit{observables} in physics. Adopting such a heavy and relatively enriched language allows us to formalize the notions of \textit{quantization} and \textit{observables} in quantum field theory as well. This document is organized to examine each task listed above in an expository manner. 
	\end{abstract}

\tableofcontents

\section{Introduction}	

Investigating higher structures in geometry and physics has undoubtedly played several important roles in conceptualizing certain problems in different branches of mathematics; ranging from algebra, number theory, topology and geometry to mathematical physics. To be more precise, it indeed leads to new models and frameworks to reformalize the existing theories and the relevant problems so that a number of new mathematical apparatuses turn out to be applicable. As  seen in the literature (see, for instance, \cite{Anel,Toen,Vezz}), one always needs alternative perspectives, new mathematical tools, and hence a different level of abstraction not only to deal with standstill problems but also to create new windows for the solutions to prospective problems. 

In this paper, we essentially concentrate on  \textit{higher/derived spaces} and geometric structures induced by higher algebraic structures. To this end, we shall rely on the language of higher category theory and derived algebraic geometry. There are, in fact, equivalent ways of formalizing derived algebraic geometry  (DAG); for instance, the one by To$\ddot{e}$n and Vezzosi \cite{Toen-Vezzosi}, and the other is Lurie's work \cite{Luriethesis}. 

DAG essential provides a new setup to deal with non-generic situations (e.g. non-transversal intersections and ``bad" quotients) in geometry. To this end, it combines higher categorical objects and homotopy theory with many tools from homological algebra. Hence, roughly speaking, it can be considered as a \emph{higher categorical/homotopy theoretical refinement of classical algebraic geometry}. In that respect, it offers a new way of organizing information for various purposes. Therefore, it has many interactions with other mathematical domains \cite{Anel,Toen}. 

DAG, for instance, provides new mathematical formalisms to capture certain aspects of intersection theory and moduli problems. We indeed closely follow this direction for the following reason:  Moduli theory in fact plays a significant role in analyzing field theories as well. Throughout the discussion we shall focus on the classical side of the story, indeed. 

The reason why one may prefer to adopt a moduli theoretic approach to study field theories is relatively simple: A \emph{classical field theory} can be described by a piece of data $(M, F_M, S, G)$, where $F_M$ denotes the \emph{space of fields} on some base manifold $M$, $\mathcal{S}$ is a \emph{smooth action functional} on $F_M$, and $G$ is a certain \emph{group encoding the symmetries} of the system \cite{Mnev}. Then the standard folklore suggests that the key information about the system is encoded in the \textit{critical locus} $ crit(\mathcal{S}) $ of $\mathcal{S}$, modulo symmetries. Therefore, the problem of interest boils down to an analysis of the properties of this moduli space. 

Of course, moduli theory has some natural questions related to ``bad quotients" or ``bad intersections" \cite{Toen}.  There are some classical techniques to deal with these sorts of problems, but nowadays some people prefer to use a relatively new technology, namely \textit{derived algebraic geometry}. 
Consequently, DAG suggests new and alternative perspectives in physics as well. In that respect, for instance, the formulation of gauge theories in the language of derived algebraic geometry provides new and fruitful insights to encode the formal geometry of the associated moduli space of the theory. Some works in this direction are \cite{Benini,Calaque,Cosv1,Cosv2,Costello Renormalizetion book,Costello Renormalisation and BV formalism paper,Gwilliam}. 

An outline of the remainder of this paper is as follows: In Sections \ref{section_setup} and \ref{structure of observables}, to make the first touch with physics and realize where derived geometry comes into play, we shall first discuss certain notions and structures in a rather intuitive manner, such as (higher) spaces of solutions to the field equations of theories and factorization algebras of observables in these theories. Afterwards, in Sections \ref{section_derived formulations} and \ref{sections_recasting some examples}, we shall investigate \textit{derived} formulations of field theories and their consequences.


\section{Classical setup and motivation} \label{section_setup}

We first recall how to define a classical field theory in Lagrangian formalism \cite{Mnev}: 

\begin{definition} \label{defn_naive defn of cft}
	A \textit{classical field theory} on a manifold $M$ consists of the following data: \begin{itemize}\itemsep=2pt
		\item [\textit{(i)}] the space $ \mathbb{F}_{M}$ of \textit{fields} of the theory defined to be the space $\Gamma(M,\mathcal{F})$ of sections of a particular \textit{sheaf} $ \mathcal{F} $ on $M$,
		\item [\textit{(ii)}] the action functional $\mathcal{S}: \mathbb{F}_{M}\longrightarrow k ~(\mathbb{R}~or~ \mathbb{C})$  
	\end{itemize}
\end{definition}  

\begin{remark} \label{rmk on critical locus}
	In order to encode the dynamics of the system, we need to study the \textit{critical locus} $ crit(\mathcal{S}) $ of $\mathcal{S}$. One can determine $ crit(\mathcal{S}) $ by employing variational techniques for the functional $\mathcal{S}$ and that leads to define $ crit(\mathcal{S}) $ to be the \emph{space of solutions to the Euler-Lagrange equations (E-L) modulo gauge equivalences}. Therefore, a classical field theory can be thought of as a \textit{study of the moduli space of solutions to the E-L equations} of the theory of interest.
\end{remark}    

\begin{definition}
	A classical field theory on a manifold $M$ is called 
	\[	\begin{cases}
		scalar & \mbox{if} \ \ \ \mathbb{F}_{M}:=C^{\infty}(M), \\
		gauge  & \mbox{if} \ \ \ \mathbb{F}_{M}:= \mathcal{A}, \\
		\sigma-model & \mbox{if} \ \ \ \mathbb{F}_{M}:= Maps(M,N).
	\end{cases} \]
	Here $\mathcal{A}$ is the space  of all $ G $-connections on a principal $ G$-bundle over $M$, namely $\mathcal{A}=\Omega^1 (M) \otimes \mathfrak{g}$, and $ Maps(M,N)$ denotes the space  of smooth maps from $M$ to $N$ for some fixed target $ N$.
\end{definition}

\begin{example}
	(\cite{KW}, Part III, Section 3.1) In accordance with the above definitions we consider the underlying theory (given as a $\sigma$-model) for a classical free particle of mass $m$ moving in $ \mathbb{R}^n $ together with a   potential energy $V:\mathbb{R}^n\rightarrow\mathbb{R}$.
	
	Let $ \mathbb{F}_{M}:=Maps(M,\mathbb{R}^n)$ with $M:= [0,1]$ (in that case $ \mathcal{F} $ is just the trivial bundle on $M$) and the action functional be given as 
	\begin{equation}
		\mathcal{S}(q):= \int_{[0,1]} \bigg(\dfrac{m||\dot{q}||^2}{2}-V(q)\bigg)~~~~for~ all~~ q:[0,1] \longrightarrow \mathbb{R}^n.
	\end{equation}
	Then the corresponding Euler-Lagrange equation becomes
	\begin{equation}
		m\ddot{q}=-grad~ V(q),
	\end{equation}
	which is indeed \textit{the Newton's equation of motion}.
	
\end{example} 
\begin{example}
	Consider a classical free particle (of unit mass) moving in a Riemannian manifold $ N $ without any potential energy. Set $ \mathbb{F}_{M}:=Maps(M,N)$ with $M:= [0,1]$. Let $f\in Maps(M,N)$ be a smooth path in $N$, and the action functional given by 
	\begin{equation}
		\mathcal{S}(f):= \frac{1}{2}\int_{[0,1]} ||\dot{f}||^2.
	\end{equation}
	This is in fact called the \textit{energy functional}  in  Riemannian  geometry. Then the corresponding Euler-Lagrange equations in a local chart $x=(x^j)_{j=1,...,dimN}$ for $N$ are given as
	\begin{equation}
		\ddot{f}^k+\Gamma^k_{i j} \dot{f}^i\dot{f}^j=0~~for ~ k=1,2,...,dimN,
	\end{equation}
	where $f^k$ denotes local component of $f$, i.e, $f^k:=x^k\circ f$, and $ \Gamma^k_{i j}:=\Gamma^k_{i j}(f(t))$ is the Christoffel symbol for each $i,j,k$.  These equations are indeed the \textit{geodesic equations in Riemannian geometry}.
	
\end{example}

\begin{example} \label{example_free scalar massive field theory}
	(\cite{Mnev}, Sec. 1.2.) Consider the theory with free scalar massive fields. Let $M$ be a Riemannian manifold and set $\mathbb{F}_{M}:=C^{\infty}(M)$. Let $\phi\in \mathbb{F}_{M}$, then we define the action functional governing the theory as 
	\begin{equation}
		\mathcal{S}(\phi):= \int_{M} \bigg(\frac{||d\phi||^2}{2}-\frac{m^2}{2}\phi^2\bigg).
	\end{equation} 
	The corresponding E-L equation in this case reads as 
	\begin{equation}
		(\Delta+m^2)\phi=0.
	\end{equation}
	
\end{example}
\begin{example} \label{example_CS theory}
	(Section 2 of \cite{Mnev}) Consider the classical $SU(2)$-Chern-Simons gauge theory on a closed, orientable 3-manifold $ X $.  The main ingredients of this structure are encoded by the theory of principal $ G $-bundles: Let $\pi: P\rightarrow X$ be a trivial principal $SU(2)$-bundle on $X$.  
	Then there exists a globally defined nowhere vanishing section $ \sigma \in \Gamma(X,P)$ commuting the diagram
	\begin{equation}
		\begin{tikzpicture}
			\matrix (m) [matrix of math nodes,row sep=3em,column sep=4em,minimum width=2em] {
				P  & P \\
				\ & X \\};
			\path[-stealth]
			(m-1-1)	edge  node [above] {$ \bullet \ SU(2) $}  (m-1-2)
			(m-1-2) edge node [right] {$\pi$} (m-2-2)
			(m-2-2) edge [bend left=40] node [left] {$\sigma$} (m-1-2);
		\end{tikzpicture}
	\end{equation} where $ \bullet \ SU(2) $ denotes the action of $SU(2)$ on the space $P$.
	
	Assume $\omega$ is a Lie algebra-valued connection 1-form on $P$. Let $ A:= \sigma^* \omega $ be its  representative, called the \textit{Yang-Mills field}. Then for this theory, the space $ \mathbb{F}_{X} $ of \textit{fields} is defined to be the space of all $\mathfrak{su}(2)$-valued connection 1-forms on $X$. That is, $\mathbb{F}_X:= \Omega^1(X,\mathfrak{su}(2))$. We denote this space simply by $\mathcal{A}$. Note that, in this case, the corresponding sheaf $ \mathcal{F} $ is the \textit{``twisted"} cotangent bundle $T^* X\otimes (X \times \mathfrak{su}(2))$, and $\mathbb{F}_X$ is the space of sections of this tensor product.
	
	Now, for this setup, the \textit{Chern-Simons action functional} $ CS: \mathcal{A} \longrightarrow S^1$ is given by
	\begin{equation}
		CS(A):=\frac{k}{4\pi} \displaystyle \int \limits_{X} \mathrm{Tr}(A\wedge \mathrm{d} A +\frac{2}{3}  A \wedge A \wedge A), ~~~~ k\in \mathbb{Z},
	\end{equation}together with the gauge group $\mathcal{G}=Map(X,SU(2))$ acting on the space $\mathcal{A}$ as follows: For all $g\in \mathcal{G}$ and $A \in \mathcal{A}$, we set \begin{equation}
		g\triangleleft A := g^{-1}\cdot A \cdot g + g^{-1}\cdot \mathrm{d} g.
	\end{equation} 
\newpage
Here, the integer $k$, called the \textit{level} of Chern-Simons theory, is related to gauge invariance properties of the quantized theory. As stressed in \cite{Mnev}, the integer $k$ plays a crucial role in the consistency of the quantized theory, which does not require single-valuedness of the action $CS$, but only that of the function $\exp(iCS)$. That is, under a general gauge transformation, $CS(A)$ can change by an integer multiple of $2\pi$; whereas, the function $\exp(iCS)$, which is the main object of interest in path integral quantization,  is a gauge invariant function on $\mathbb{F}_X.$ 

The argument above essentially relies on the description of $CS$ in terms of characteristic classes and geometric invariants. In that respect, the action of the Chern-Simons theory is in fact proportional to the integral of the so-called \emph{CS 3-form} $$\mathrm{Tr}(A\wedge \mathrm{d} A +\frac{2}{3}  A \wedge A \wedge A),$$ 
	which is indeed obtained from a particular characteristic class, namely the \textit{Chern form $P_4(F_A)$}. Here $F_{A}=\mathrm{d} A+A \wedge A$ is the $ \mathfrak{su}(2)$-valued curvature 2-form on $X$ associated to $A\in\Omega^1(X,\mathfrak{su}(2)).$ More details on the construction of $CS$ $(2n-1)$-forms by means of Chern forms $P_{2n}(F_A)$ and their gauge invariance properties can be found in \cite{Zanelli}.
	
	
	From the data above, the corresponding Euler-Lagrange equation in this case turns out to be $$F_{A}=0.$$  
	Note that the moduli space $ \mathcal{M}_{flat}$ of flat connections, i.e. $A\in \mathcal{A}$ with $F_{A}=0$ modulo gauge transformations, emerges in many other mathematical domains; such as topological quantum field theory, low-dimensional quantum invariants or (infinite dimensional) Morse theory.
\end{example}

\paragraph{Towards derived formulations.}Definition \ref{defn_naive defn of cft} above can be re-stated by using the language of \textit{deformation theory and derived geometry.} For the sketch of essential ideas, we refer to Chapter 3 of \cite{Cosv2}. 
The standard folklore suggests that for a classical field theory, one can make a reasonable measurement \textit{only on} those fields which are the solutions to the Euler-Lagrange equations of the given action functional. Therefore, \textit{measurements or observables} are those functions defined on the space $ \mathcal{EL} $ of solutions to the Euler-Lagrange equations, and hence by adopting the Lagrangian formalism, a classical field theory can be thought of as a study of the \textit{critical locus} of the action functional as indicated in Remark \ref{rmk on critical locus}.

However, in order to circumvent certain issues with the critical locus, which will be elaborated later in the text, we work with the \textit{derived moduli space of solutions} instead of the naïve one. Furthermore, we usually intend to capture the \textit{perturbative} behavior of the theory at the same time. Hence, this derived moduli space will be defined as the so-called \textit{formal moduli problem}. In other words, a classical field theory can be realized as a \textit{formal moduli problem} (in the sense of \cite{Lurie}) cut out by a system of PDEs determined by the corresponding Euler-Lagrange equations. 

Note that instead of a na\"{\i}ve moduli we consider it as a \textit{derived} moduli problem (the reason will be discussed below) and we would like to understand the \textit{local} classical observables $Obs^{\textit{cl}}(U)$ for all open subset $ U \subseteq M $ as well. Therefore, the assignment 
\begin{equation}
	U\xrightarrow{\mathcal{EL}}\mathcal{EL}(U)
\end{equation}
can be realized as a \textit{sheaf} of  derived spaces of solutions to the E-L equations (i.e. a sheaf of certain \textit{derived stacks}), and hence the space $ Obs^{\textit{cl}}(U) $ of classical observables on $U$ is defined as 
\begin{equation}
	Obs^{\textit{cl}}(U):=\mathcal{O}_{\mathcal{EL}(U)},
\end{equation} 
where $ \mathcal{O}_{\mathcal{EL}(U)} $ denotes the \emph{algebra of functions} on the formal moduli space $ \mathcal{EL}(U) $. 

In the language of derived schemes/stacks, derived objects are locally modeled on commutative differential graded algebras (\textit{cdgas}), and hence in our case $ Obs^{\textit{cl}}(U) $ can naturally be realized as a certain cdga.  Moreover, the space $ Obs^{\textit{cl}}(U) $ is the dual space of $ \mathcal{EL}(U) $ for each open subset $ U \subseteq M $, and hence the assignment 
\begin{equation}
	U\xrightarrow{Obs^{\textit{cl}}} Obs^{\textit{cl}}(U)
\end{equation} 
gives rise to a certain\textit{ co-sheaf,} which will be discussed in a rather succinct and na\"{\i}ve way below.
\newpage
\section{Factorization algebras and the structure of observables} \label{structure of observables} Main motivations in \cite{Cosv1} and \cite{Cosv2} to study \textit{factorization algebras} associated to a \textit{perturbative} QFT is to generalize the deformation quantization approach to quantum mechanics developed by Kontsevich \cite{Kontsevich}. In other words, deformation quantization essentially encodes the nature of observables in \textit{one}-dimensional quantum field theories, and the \textit{factorization algebra formalism} provides an $ n $-dimensional generalization of this approach. 

To be more precise, recall that observables in classical mechanics and those in corresponding quantum mechanical system can be described in the following way: Let $(M,\omega)$ be a symplectic manifold (a \textit{phase space}), then we define  the \textit{space  $A^{\textit{cl}}$ of classical observables on $M$} to be the space $C^{\infty}(M)$ of smooth functions on $M$. Hence, $A^{\textit{cl}}$ forms a \textit{Poisson algebra} with respect to the Poisson bracket $\{ \cdot, \cdot \}$ on $C^{\infty}(M)$ given by 
\begin{equation} \label{defn of poisson bracket}
	\{f,g \} := -w(X_f,X_g)=X_f (g) \ \ for \ all  \ f,g \in C^{\infty}(M),
\end{equation} where  $X_f$ is the \textit{Hamiltonian vector field associated to} $f$ defined implicitly as \begin{equation} \label{defn of Hamiltonian vector field}
	\imath_{ X_f} \omega=df.
\end{equation} Here, $ \imath_{ X_f} \omega $ denotes the \textit{contraction} of a 2-form $\omega$ with the vector field $X_f$ in the sense that 
\begin{equation}
	\imath_{ X_f}  \omega \ (\cdot):= \omega(X_f, \cdot).
\end{equation}Employing canonical/geometric quantization formalism \cite{Blau}, the process of quantization boils down to  \emph{studying representation theory of (a certain subalgebra $\mathcal{A}$ of) classical observables} in the sense that one can construct the \textit{quantum Hilbert space $\mathcal{H}$} and a particular homomorphism \footnote{The homomorphism $\mathcal{Q}$ indeed becomes a Lie algebra homomorphism after a small modification: Recall that \textit{a Lie algebra homomorphism} $\beta:\mathfrak{g} \rightarrow \mathfrak{h}$ is a linear map of vector spaces such that $\beta([X,Y]_{\mathfrak{g}})=[\beta(X),\beta(Y)]_{\mathfrak{h}}$. Keep in mind that, one can easily suppress the constant ``-$i\hbar$" in Equation \ref{quantum cond.} into the definition of $\mathcal{Q}$ such that the quantum condition in Equation \ref{quantum cond.} becomes the usual compatibility condition for being a Lie algebra homomorphism.} \begin{equation}
	\mathcal{Q}: \mathcal{A}\subset \big(C^{\infty}(M),\{ \cdot, \cdot \}\big)\longrightarrow \big(End(\mathcal{H}),[\cdot , \cdot ]\big)
\end{equation}
together with the \textit{Dirac's quantum condition}: $ \forall $ $f,g\in \mathcal{A}
$ we have
\begin{equation} \label{quantum cond.}
	[\mathcal{Q}(f), \mathcal{Q}(g)]=-i\hbar \mathcal{Q}\big(\{f ,g \}\big),
\end{equation} where $[\cdot, \cdot]$ denotes the usual commutator on $End(\mathcal{H})$.

In accordance with the above setup, while the classical observables form a \textit{Poisson algebra}, the space $A^{\textit{q}}$ of quantum observables forms an \textit{associative algebra} which is related to classical one by the \textit{quantum condition} given in Equation \ref{quantum cond.}. Deformation quantization, in fact, serves as a mathematical treatment that captures this correspondence, i.e. \textit{the deformation of a commutative structure to a non-commutative one}, for general Poisson manifolds. For details, see \cite{Kontsevich}.

Factorization algebras, on the other hand, are algebro-geometric objects which are manifestly described in the language analogous to that of \textit{(co-)sheaves}. For a complete discussion, see Ch. 3 of \cite{Cosv1}, or \cite{Ginot}.

\begin{definition} \label{defn of prefactorization algebra}
	A \textit{prefactorization algebra $\mathcal{F}$} on a manifold $M$ consists of the following data: \begin{itemize}
		\item For each open subset $ U \subseteq M $, a cochain complex $\mathcal{F}(U)$. 
		\item For each open subsets $ U \subseteq V $ of $M$, a cochain map $ \imath_{U;V} :  \mathcal{F}(U)\longrightarrow\mathcal{F}(V)$.
		\item For any finite collection $U_1, ..., U_n$ of pairwise disjoint open subsets of $V\subseteq M$, $V$ open in $M$, there is a morphism  \begin{equation}
			\imath_{U_1,...U_n;V} :  \mathcal{F}(U_1) \otimes \cdot \cdot \cdot \otimes \mathcal{F}(U_n)\longrightarrow\mathcal{F}(V)
		\end{equation} together with \textit{certain compatibility conditions:}  
		\vspace{5pt}
		
		\begin{enumerate}
			\item[\textit{i.}] The \textit {invariance} under the action of symmetric group $S_n$ permuting the ordering of the collection $U_1, ..., U_n$ in the sense that \begin{equation}
				\imath_{U_1,...U_n;V} = \imath_{U_{\sigma(1)},...U_{\sigma(n)};V} \ for \ any \ \sigma \in S_n.
			\end{equation} That is, the morphism $\imath_{U_1,...U_n;V}$ is independent of the ordering of open subsets $U_1, ..., U_n$, but it depends only on the family $\{U_i\}$.
			\item[\textit{ii.}]  The \textit{associativity} condition in the sense that if $U_{i1} \ \amalg \cdot \cdot \cdot \amalg \ U_{in_i} \subset V_i$ and $V_{1} \ \amalg \cdot \cdot \cdot \amalg \ V_{k} \subset W$ where $U_{ij} (resp. \ V_i)$ are pairwise disjoint open subsets of $V_i$ (resp. $W$) with $W$ open in $M$, then the following diagram commutes.
			\begin{equation}
				\begin{tikzpicture}
					\matrix (m) [matrix of math nodes,row sep=3em,column sep=4em,minimum width=2em] {
						\bigotimes\limits_{i=1}^{k}	\bigotimes\limits_{j=1}^{n_i} \ \mathcal{F}(U_{ij})   & \bigotimes\limits_{i=1}^{k} \mathcal{F}(V_i) \\
						& \mathcal{F}(W) \\};
					\path[-stealth]
					(m-1-1) edge node [left] {} (m-2-2)
					edge  node [below] {} (m-1-2)
					(m-1-2) edge node [right] {$ $} (m-2-2);
				\end{tikzpicture}
			\end{equation}
		\end{enumerate} 
	\end{itemize}
	
\end{definition}  
With this definition in hand, a prefactorization algebra behaves like a co-presheaf except the fact that we use tensor product instead of a direct sum of cochain complexes. Furthermore, we can define  a \textit{factorization algebra} once we impose a certain local-to-global condition on a prefactorization algebra analogous to the one imposed on presheaves. For details, see Ch. 4 \& Ch. 8 of \cite{Ginot}. 

Factorization algebras in fact serve as $ n $-dimensional counterparts to those objects realized in deformation quantization formalism. In particular, one recovers the observables in classical/quantum mechanics when restricts to the case $n=1$ (for a readable discussion see Ch. 1 of \cite{Cosv1}). For instance, in a particular gauge theory, holonomy observables, namely \textit{Wilson line operators}, can be formalized in terms of such objects. For details, we refer to Part 3, Ch. 8.2 of \cite{Cosv1}. These are, in fact, the ones that are related to Witten's Knot invariants. They actually arise from the analysis of certain partition functions in three-dimensional Chern-Simons theory. In this approach to perturbative quantum field theories, quantum observables in these types of theories form a factorization algebra (e.g. Lemma 2.3.2 in \cite{Cosv1}). 

Last, but not least, it also turns out that a factorization algebra of quantum observables is related to a (commutative) factorization algebra of associated classical observables in the following sense: 
\begin{theorem}\label{thm_weekquan}
	(Weak quantization Theorem \cite{Cosv2} Sec. 1.3): For a classical field theory and a choice of BV quantization,\begin{enumerate}
		\item The space $Obs^{\textit{q}}$ of quantum observables forms a factorization algebra over the ring $ \mathbb{R}\big[[\hbar]\big] $.
		\item $Obs^{\textit{cl}} \cong Obs^{\textit{q}}~mod ~\hbar$ as a homotopy equivalences where $Obs^{\textit{cl}}$ denotes the associated factorization algebra of classical observables.
	\end{enumerate}
\end{theorem}

\begin{remark} As argued in Section 2 of \cite{Cos2}, for any open subset $U$,  $Obs^{q}$ is a certain cohain complex  of $\mathbb{R}\big[[\hbar]\big]$-modules, flat over $\mathbb{R}\big[[\hbar]\big]$ such that these cochain complexes are, informally speaking, obtained from spaces of smooth functions and distributions on the spacetime manifold under consideration. Furthermore, the second statement of Theorem \ref{thm_weekquan} implies that the factorization algebra of quantum observables deforms that of classical observables in a way that classical observables can be obtained from the quantum ones by taking modulo $\hbar$ \cite{Cos2}. Notice that  $Obs^{\textit{cl}}(U) \cong Obs^{\textit{q}}(U)~mod ~\hbar$ is just the \textit{derived} space of functions on the \textit{derived} moduli space of solutions to the E-L equations over $U$. 
\end{remark}

Note that the  theorem above is just a part of the story, and it is indeed \textit{weak} in the sense that it is \textit{not} able to capture the data related to Poisson structures on the space of observables. To provide a correct $ n $-dimensional analogue of deformation quantization approach, we need to refine the notion of a classical field theory in such a way that the richness of this new set-up become visible. This is where \textit{derived algebraic geometry} comes into play. As we discussed above, the space of classical observables forms a (commutative) factorization algebra $Obs^{cl}$, where $Obs^{cl}(U)$ is just the \textit{derived} space of functions on the \textit{derived} moduli space of solutions to the E-L equations over $U$. In that respect, to encode the desired Poisson structures, one may employ the theory of derived Possion geometry \cite{CPTVV, Melani, Safranov}. 

The structure of a factorization algebra allows us to employ certain \textit{cohomological methods} encoding that of observables in a given theory in the following sense (cf. \cite{Cosv1} Ch. 1, Sec. 3): Factorization algebra $ Obs^{\textit{cl}} $ of observables can be realized as a particular assignment analogous to a co-sheaf of \textit{cochain complexes} as mentioned above. That is, for each open $U\subseteq M$, $ Obs^{\textit{cl}} (U)$ has a $ \mathbb{Z} $-graded structure
$$ Obs^{\textit{cl}}(U)= \bigoplus_{i\in \mathbb{Z}} Obs^{\textit{cl}}_{i}(U) 
$$
together with suitable connecting homomorphisms $\mathrm{d}_{i}:Obs^{\textit{cl}}_{i}(U)\rightarrow Obs^{\textit{cl}}_{i+1}(U) $ for each $i$. The corresponding cohomology groups $H^{i} (Obs^{\textit{cl}} (U))$ encode the structure of observables as follows:
\begin{itemize}
	\item \textit{``Physically meaningful"} observables are the closed ones with cohomological degree 0, i.e., $ \mathcal{O} \in Obs^{\textit{cl}}_{0} (U) $ with $ \mathrm{d}_{0} \mathcal{O}=0. $ (and hence $ [\mathcal{O}] \in H^{0} (Obs^{\textit{cl}} (U) $).)
	\item $ H^{1} (Obs^{\textit{cl}} (U)) $ contains \textit{anomalies}, i.e., obstructions for classical observables to be lifted to the quantum level. In a gauge theory, for instance, there exists certain classical observables respecting gauge symmetries such that they do \textit{not} admit any lift to quantum observables respecting gauge symmetries. This behaviour is indeed encoded by a non-zero element in $ H^{1} (Obs^{\textit{cl}} (U)) $  
	\item $ H^{n} (Obs^{\textit{cl}} (U)) $ with $ n<0 $ can be interpreted as symmetries, higher symmetries of observables etc. via higher categorical arguments.
	\item $  H^{i} (Obs^{\textit{cl}} (U)) $ with $ n>1 $ has no clear physical interpretation.

\end{itemize} 

\section{Derived formulation of field theories} \label{section_derived formulations}
Using the \textit{derived} interpretation of a classical field theory
, one can employ a number of mathematical techniques and notions naturally appear in derived algebraic geometry. For instance, we may consider a classical field theory as a study of  the \textit{derived critical locus} \cite{Joyce2, Vezz} of the action functional since it can be considered as a formal moduli problem in the sense indicated above. Indeed, passing to the \textit{derived} moduli space of solutions corresponds to Batalin-Vilkovisky formalism for  classical field theories \cite{Calaque}. 

In derived algebraic geometry, any formal moduli problem arising as the derived critical locus admits a $ -1 $-shifted symplectic structure \cite{Calaque,PTVV,Vezz}. This observation is crucial and it ensures the existence of a symplectic structure on the space $ Obs^{\textit{cl}} $ of classical observables. In the language of derived algebraic geometry, therefore, we have the following definition (Appendix of \cite{Cos2}, Chapter 3 of \cite{Cosv2}):
\begin{definition} \label{main defn}
	Given  $(M, F_M, S, G)$,  where $F_M$ denotes the space of fields on some base manifold $M$, $\mathcal{S}$ is a smooth action functional on $F_M$, and $G$ is a certain symmetry group acting on $F_M$, then we define a \emph{\textbf{(perturbative) classical field theory}} on $M$ to be a sheaf of \textit{derived stacks} (of the \textit{derived critical locus} $ dcrit(\mathcal{S})$) with a symplectic form of degree $ -1 $.
\end{definition} 

\begin{remark}
	Definition \ref{main defn} indicates that a classical field theory can be ideally described as a \emph{sheaf of derived stacks of solutions to the EL-equations} with an appropriate \emph{derived} symplectic/Poisson structure. But, when we require to encode the perturbative behavior, it is better to use \emph{formal moduli problems} instead of derived stacks. Here, formal moduli problems are indeed particular derived stacks with more specific local algebraic structures that are used as suitable models to capture the perturbative nature of the problem under consideration.	
\end{remark}
 We now intend to \textit{unpackage}  Definition \ref{main defn} in an intuitive way. The remaining part of the current section will be devoted to that purpose.
\subsection{Why does the term ``derived" emerge?}  To explain how ``derived" techniques enter the picture, we  first need to discuss the na\"{\i}ve or underived realization of a classical field theory in the language of \textit{intersection theory}. 

Let $F_M$ denote the space  of fields on a base manifold $M$ as in Definition \ref{main defn}. For simplicity, we assume $F_M$ is a finite dimensional manifold. As indicated in the Remark \ref{rmk on critical locus}, a classical field theory can be encoded by the critical locus $ crit(\mathcal{S}) \subset F_M$ of the action functional $\mathcal{S}$ on $F_M$. However, computing certain path integrals perturbatively around classical solutions to the E-L equations is usually problematic if the critical points are \textit{degenerate}. To avoid such problems, one can employ a certain  trick the so-called \textit{Batalin-Vilkovisky formalism} which, roughly speaking, consists of adding certain fields, such as ghosts, anti-fieds etc..., to the functional \cite{Calaque}. This problem, on the other hand, can be formulated in the language of intersection theory as follows. 

Note that $ crit(\mathcal{S}) $  is defined to be the  \emph{intersection of the graph $\Gamma(\mathrm{d}S)$ of $\mathrm{d}S$ with the zero-section of the cotangent bundle $T^*F_M$ inside $T^*F_M$} (Chapter 5 of \cite{Cosv2}). That is,
\begin{equation}
	crit(\mathcal{S}) :=  \Gamma(\mathrm{d}S) \cap F_M, \text{ where } \mathrm{d}S \in \Omega^1(F_M).
\end{equation}

\noindent By adopting an algebro-geometric language, $ crit(\mathcal{S}) $ can be described in terms of the \textit{fibered product}  $ F_M \times_{T^*F_M} F_M $ such that the following diagram commutes:
\begin{equation}
	\begin{tikzpicture}
		\matrix (m) [matrix of math nodes,row sep=3em,column sep=4em,minimum width=2em] {
			crit(\mathcal{S}):= F_M \times_{T^*F_M} F_M  & F_M \\
			F_M & T^*F_M \\};
		\path[-stealth]
		(m-1-1) edge node [left] {} (m-2-1)
		edge  node [below] {} (m-1-2)
		(m-2-1.east|-m-2-2) edge node [below] {} node [above] {$\mathrm{d}\mathcal{S}$} (m-2-2)
		(m-1-2) edge node [right] {$0$} (m-2-2);
	\end{tikzpicture}
\end{equation}

\noindent Even if $F_M$ is a smooth manifold, for instance, the intersection $F_M \times_{T^*F_M} F_M $ may be highly problematic, and hence an object $ \big(crit(\mathcal{S}),\mathcal{O}_{crit(\mathcal{S})} \big)$ generically fails to live in the same category; i.e. intersection would not define a manifold at all (e.g. non-transverse intersection of  two submanifold is not a submanifold in general).  In fact, problems arising from the degeneracy of critical points correspond to problematic intersection in the above sense \cite{Calaque}. 

If we employ, however, a derived setup and introduce  derived geometric counterparts of the objects under consideration, such as \textit{derived schemes or dg-schemes}, then we can circumvent the non-existence problem for fibered products. 
For simplicity, one can work with a dg-scheme $ (X,\mathcal{O}_{X}) $ for which the structure sheaf $ \mathcal{O}_{X}$ is a sheaf of commutative differential graded algebras \cite{Arinkin}. 

In general, a dg-scheme is a more structured object than a derived scheme. Indeed, there exists a forgetful functor from dg-schemes to derived schemes, which is neither full nor faithful nor essentially surjective \cite{Arinkin,  Toen}. But, in the case of \textit{embeddable derived schemes} (derived schemes that can be embedded in a scheme), dg-schemes can provide  models for such derived schemes. For details, see Section 3.1 of \cite{Toen}. Furthermore, one can show that the $\infty$-category of derived schemes (resp. the 2-category dg-schemes) admits  \textit{fibered products}. Therefore, this leads to the following motivation behind the use of \textit{``derived"} formulation in Definition \ref{main defn}: 

\begin{quote}
	\textit{One has to provide a new framework in which, after performing operations like a fibered product or taking quotients, the resulting objects must stay in the same place with the original objects of interest.}
\end{quote}
This requires to re-organize the \textit{local model} for the intersection of ringed spaces as follows: As in \cite{Arinkin} or Chapter 1 of \cite{Luriethesis}, instead of a na\"{\i}ve intersection determined algebraically by 
\begin{equation}
	\mathcal{O}_{crit(\mathcal{S})}:=\mathcal{O}_{G(dS)} \otimes_{\mathcal{O}_{T^*M}} \mathcal{O}_{M},
\end{equation} 
we introduce the derived version:
\begin{equation}
	\mathcal{O}_{dcrit(\mathcal{S})}:=\mathcal{O}_{G(dS)} \otimes_{\mathcal{O}_{T^*M}}^{\mathbb{L}} \mathcal{O}_{M}.
\end{equation} 
where $\cdot \otimes_{\mathcal{O}_{T^*M}}^{\mathbb{L}} \cdot $ denotes the \textit{derived tensor product}.
\vspace{5pt}

\noindent\textit{\textbf{A digression on the definition of $\cdot \otimes_{\mathcal{O}_{T^*M}}^{\mathbb{L}} \cdot $.}} To understand the basic properties of this item, we follow Chapter 0 of \cite{LurieSAG}: Let $R$ be a commutative ring, $B$ a $R$-module. Then the \textit{derived tensor product}   $ \cdot \otimes_{R}^{\mathbb{L}} B $ arises from the construction of the \textit{left-derived functor} associated to the right-exact functor 
\begin{equation}
	\cdot \otimes_{R} B: Mod_{R} \rightarrow Mod_{R}.
\end{equation}
Let $A$, $B$ be two commutative algebras over $R$. Then the definition of $ A \otimes_{R}^{\mathbb{L}} B $ naturally appears in the construction of the $i^{th}$ Tor groups $Tor_{i}^{R}(A,B)$ given by the $i^{th}$ homology of the \textit{\textbf{tensor product complex}}  $(P_{\bullet}\otimes_{R} B,\mathrm{d'})$:
\begin{equation}
	\cdots\longrightarrow P_2 \otimes_{R} B \longrightarrow P_1 \otimes_{R} B\xrightarrow{\mathrm{d'}} P_0 \otimes_{R} B \longrightarrow 0,
\end{equation}
where $P_{\bullet}$ is a projective resolution of $A$ equipped with a differential $\mathrm{d}$ such that $(P_{\bullet},\mathrm{d})$ becomes a commutative dg-algebra over $ R $ and $\mathrm{d'}=\mathrm{d}\otimes_{R} id_{B}$.
Since $B$ is a commutative $R$-algebra, the tensor product complex inherits the structure of a commutative dg-algebra over $R$ as well, and we denote this tensor product complex by $ A \otimes_{R}^{\mathbb{L}} B $. That is, we set
\begin{equation}
	A \otimes_{R}^{\mathbb{L}} B := (P_{\bullet}\otimes_{R} B,\mathrm{d'}).
\end{equation}

\begin{remark}
	The resulting commutative dg-algebra $ A \otimes_{R}^{\mathbb{L}} B $ is independent of the choice of $(P_{\bullet}\otimes_{R} B,\mathrm{d'})$ up to a quasi-isomorphism. \\
	\noindent \textit{The end of digression.}
\end{remark}

If we go back to the local model discussion for the derived tensor products, then for each open subset $ U \subseteq M $ we have
\begin{equation}
	\mathcal{O}_{dcrit(\mathcal{S})}(U):=\mathcal{O}_{G(dS)}(U) \otimes_{\mathcal{O}_{T^*M}(U)}^{\mathbb{L}} \mathcal{O}_{M}(U),
\end{equation} where RHS corresponds to the tensor product complex of dg-algebra as above. 

Together with this local model, $ \big(dcrit(\mathcal{S}),\mathcal{O}_{dcrit(\mathcal{S})} \big)$ becomes a dg-scheme with its structure sheaf $ \mathcal{O}_{dcrit(\mathcal{S})} $ being the sheaf of commutative dg $k$-algebras such that
$ \mathcal{O}_{dcrit(\mathcal{S})} $ can be manifestly given as \textit{a Koszul resolution} of $  \mathcal{O}_{M}$ as a module over $\mathcal{O}_{T^*M}$ \cite{Vezz}:
\begin{equation}
	\mathcal{O}_{dcrit(\mathcal{S})}: ~~ \cdot \cdot \cdot \longrightarrow \Gamma(M,\wedge^2TM) \longrightarrow \Gamma(M,TM)\xrightarrow{\imath_ {\mathrm{d}\mathcal{S}}} \mathcal{O}_{M} \longrightarrow 0,
\end{equation}
where $ \Gamma(M,\wedge^iTM) $ is the space of \textit{polyvector fields of degree} $i$ \textit{(or $i$-vector fields)} and $\imath_ {\mathrm{d}\mathcal{S}} $ denotes the contraction with $\mathrm{d}\mathcal{S}$ in the sense that for any 1-vector field $X\in \Gamma(M,TM)$ we define \begin{equation}
	\imath_ {\mathrm{d}\mathcal{S}}(X):= \mathrm{d}\mathcal{S}(X)=X\mathcal{S}.
\end{equation} Then, extending to $i$-vector fields by linearity, we set    \begin{equation}
	\mathcal{O}_{dcrit(\mathcal{S})} := \bigg(\bigoplus_{i\in \mathbb{Z}_{\leq0}} \Gamma(M,\wedge^iTM),\ \imath_ {\mathrm{d}\mathcal{S}} ~\bigg).
\end{equation} 
\begin{remark}
	\begin{enumerate}
		
		\item $ \big(dcrit(\mathcal{S}),\mathcal{O}_{dcrit(\mathcal{S})} \big)$ admits further derived structure; namely, a symplectic form of cohomological degree $-1$ (see \cite{PTVV}, Corollary 2.11). The description of this structure, however, is beyond the scope of the current discussion. For the construction, we refer to \cite{PTVV}. You may also see \cite{Brav} or \cite{Joyce} for an accessible presentation of PTVV's shifted symplectic geometry.
	
		\item The existence of  derived geometric structures will be crucial when we discuss the notion of quantization for $n$-dimensional classical field theories. Indeed, this new structure is really what we need, and it leads to $n$-dimensional generalization of what we have already had in the case of quantization of classical mechanics. Recall that observables in classical mechanics with a phase space $(X,\omega)$ form a \textit{Poisson algebra} with respect to the Poisson bracket $\{ \cdot , \cdot \}$ on $C^{\infty}(X)$, which is obtained from the symplectic structure $\omega$. With the same spirit but in a more non-trivial way, a derived/shifted symplecic form on a derived moduli space in Definition \ref{main defn} will induce a derived Poisson structure \cite{CPTVV,Melani, Safranov}. 
		\item In the derived setup, to encode local models for these higher structures, there are also Darboux-like theorems and the corresponding neighborhood theorems for derived spaces with shifted symplectic/Poisson structures. For details, we refer to \cite{Brav,CPTVV,JS}.
		
	\end{enumerate} 
\end{remark}

\subsection{Why does ``stacky" language come in?} Studying the critical locus of an action functional $\mathcal{S}$ is just one part of the story, and we already observed that in order to avoid the degenerate critical points one requires to introduce the notion of derived intersection and hence the \emph{derived critical locus}, which is well-behaved than the usual one.  For details, we again refer to \cite{Calaque,Vezz,Vezz2}. 

Other part of the story is related to moduli nature of the problem. Indeed, one requires to quotient out by symmetries while studying the solution space of the E-L equations, but the quotient space might be highly problematic as well. For instance, the action of a gauge group $\mathcal{G}$ on a manifold $X$ may not be free, and hence the resulting quotient $X/\mathcal{G}$ would not be a manifold, but it can be realized as an \textit{orbifold} $[X//\mathcal{G}]$, which is indeed a particular \textit{stack}, given by the\textit{ orbifold quotient}.
\paragraph{Moduli problems and stacks.} A \textit{moduli problem} is a problem of constructing a classifying space (or a moduli space $\mathcal{M}$) for certain geometric objects (such as manifolds, algebraic varieties, vector bundles etc...) up to their intrinsic symmetries. In moduli theory, providing the so-called ``fine" moduli space $\mathcal{M}$ for the moduli problem of interest is essential. For an accessible overview, we refer to \cite{BenZ}. A relatively complete treatment can be found in Section 1.2 of \cite{Neumann}.

Let $\mathcal{C}$ be a category, we now outline the wish-list for a ``fine" moduli space $\mathcal{M} \in Ob(\mathcal{C})$: 
\begin{enumerate}
	\item  $\mathcal{M} \in Ob(\mathcal{C})$ is supposed to serve as a \textit{parameter space} in a sense that there must be a one-to-one correspondence between the points of  $\mathcal{M}$ and the \textit{set} of isomorphism classes of objects to be classified:
	\begin{equation}
		\{points \ of \ \mathcal{M}\} \leftrightarrow \{isomorphism \ classes \}
	\end{equation}
	\item The existence of a \emph{universal classifying object},  say $\mathcal{T}$, through which all other objects parametrized by $\mathcal{M}$ can  be reconstructed. This, in fact, makes the moduli space $\mathcal{M}$ even more sensitive to the behavior of ``families" of objects on any base object $B$. Roughly speaking, this is manifested by a certain representative morphism $B\rightarrow \mathcal{M}$. 
\end{enumerate}

More precisely, we define a \emph{family over a base scheme $B$} to be a scheme $X$ relative to $B$; i.e., a morphism of schemes   \begin{equation}
	\pi: X \longrightarrow B.
\end{equation}  Furthermore, for any closed point $b\in B$, a \emph{fiber $X_b$ of $\pi$ over $b$}  is defined as the fibered product $X_b=X \times_{B} \{b \}.$

Now, if $\mathcal{M}$ is the ``fine" moduli space for the moduli problem of interest and $\mathcal{T}$ is a universal family, then for any  family $\pi: X \longrightarrow B$,  there exists a unique morphism $f:B \rightarrow \mathcal{M}$ such that one has the following fibered product diagram: 
\begin{equation}
	\begin{tikzpicture}
		\matrix (m) [matrix of math nodes,row sep=3em,column sep=4em,minimum width=2em] {
			X  & \mathcal{T} \\
			B & \mathcal{M}\\};
		\path[-stealth]
		(m-1-1) edge  node [left] {{\small $ \pi $}} (m-2-1)
		(m-1-1.east|-m-1-2) edge  node [above] {{\small $  $}} (m-1-2)
		(m-2-1.east|-m-2-2) edge  node [below] {} node [below] {{\small $ f $}} (m-2-2)
		(m-1-2) edge  node [right] {{\small $  $}} (m-2-2);
	\end{tikzpicture}
\end{equation}
That is, the family $X$ can be uniquely obtained by pulling back the universal family $\mathcal{T}$ along the morphism $f.$ 

In the language of category theory, on the other hand, a \textit{moduli problem} can be formalized as a certain functor \begin{equation}
	\mathcal{F}: \ \mathcal{C}^{op} \longrightarrow \ Sets,
\end{equation}
which is called a \textit{\textbf{moduli functor}}, where $ \mathcal{C}^{op} $ is the \textit{opposite} category of the category $\mathcal{C}$, and $ Sets $ is the category of sets. In order to make the argument more transparent, we take $\mathcal{C}$ to be the category $Sch$ of schemes. Note that 
for each scheme $ U \in Sch $, $\mathcal{F}(U)$ is the \textit{set} of isomorphism classes parametrized by $U$, and for each morphism $f:U \rightarrow V$ of schemes, we have a morphism $\mathcal{F}(f): \mathcal{F}(V) \rightarrow \mathcal{F}(U)$ of sets. 

Together with the above formalism, the existence of a fine moduli space in fact corresponds to the \textit{representability} of the moduli functor $\mathcal{F}$ in the sense that 
\begin{equation}
	\mathcal{F} = Hom_{Sch} (\cdot, \mathcal{M}) \ \ for \ some \ \mathcal{M}\in Sch.
\end{equation}If this is the case, then we say that \textit{$\mathcal{F}$ is represented by $\mathcal{M}$}.
\vspace{5pt}

In many cases, however, a moduli functor is \textit{not} representable in the category $Sch$ of schemes. This usually happens for moduli problems involving non-trivial automorphism groups. This is essentially where  \textit{stacks} come into play. Stacks are categorical generalization of sheaves, and they are of particular interest in moduli theory to deal with the aforementioned representability problem \cite{Toen}. 

A stack can also be thought of as a first instance such that the ordinary notion of a category  \textit{no longer} suffices to define such an object. To make sense of this new object in a well-established manner and enjoy the richness of this \textit{new structure}, we need to introduce a higher categorical notion, namely a \textit{2-category}.  

The theory of stacks, therefore, employs higher categorical techniques and notions in a way that it provides a mathematical treatment for the representability problem by re-defining the moduli functor as a stack. For an accessible treatment, see Section 2.1 of \cite{Neumann} or Section 1 of \cite{Toen}. 

In brief, a \textit{stack} is a particular groupoid-valued pseudo-functor with local-to-global properties,
\begin{equation}
	\mathcal{X}: \mathcal{C}^{op} \longrightarrow \ Grpds,
\end{equation} 
where \textbf{$Grpds$} denotes the 2-category of groupoids with objects being categories $\mathcal{C}$ in which all morphisms are isomorphism (these sorts of categories are called \textit{groupoids}). Here, 1-morphisms are functors $\mathcal{F}: \mathcal{C} \rightarrow \mathcal{D}$ between groupoids, and 2-morphisms are  natural transformations $\psi: \mathcal{F} \Rightarrow \mathcal{F'}$ between two functors. 

\begin{remark}
	\label{rem5.1}In order to make sense of local-to-global (or ``glueing") type arguments, one requires to introduce an appropriate notion of \textit{topology on a category $\mathcal{C}$}. Such a structure is called  the \textit{Grothendieck topology} $\tau$. Furthermore, a category $\mathcal{C}$ equipped with a Grothendieck topology $\tau$ is called \textit{site.} Note that if we have a site $\mathcal{C}$, then we can define a \textit{sheaf on $\mathcal{C}$} in a well-established manner as well. This essentially leads to the \textit{functor of points}-approach for defining a scheme $X$ in the following sense: Given a scheme $X$,  one can define a sheaf (on the category $Sch$ of schemes) by using the Yoneda functor $ Hom_{Sch} (\cdot, X) $ as
	\begin{equation}
		\underline{X}: \ Sch^{op} \longrightarrow \ Sets,
	\end{equation} where $\underline{X}:=Hom_{Sch} (\cdot, X)$. This is indeed a sheaf by the theorem of Grothendieck (cf. Theorem 1.29 in \cite{Neumann}).
\end{remark}

\begin{remark}
	Any 1-category (i.e. a standard category) can be realized as a 2-category in which there exist \textit{no} non-trivial higher structures, i.e. 2-morphisms in a 1-category are just identities.
	
\end{remark}

\begin{remark}
	By using the 2-categorical version of Yoneda's lemma, one can show that the moduli functor $\mathcal{X}: \mathcal{C}^{op} \longrightarrow Grpds$, \textit{a prestack over a category $\mathcal{C}$ (with the descent property)}, turns out to be representable in the 2-category $Stks$ of stacks. For the precise statement and the construction, see Section 2.1 of \cite{Neumann}. Here, objects in $Stks$ are in fact presheaves of groupoids possessing the descent/local-to-global property w.r.t. the underlying site structure. For more details, see \cite{Toen-Vezzosi}.
\end{remark}

To sum up, as in the case of derived intersections, informally speaking, we enrich the category under consideration with certain non-trivial higher structures in a way that a moduli problem would become representable in this enriched-version even if this was not the case in the first place. The price we have to pay is to adopt higher categorical dictionary leading to the change in the level of abstraction, where objects under consideration become rather counter-intuitive. Indeed, stacks and 2-categories serve as  motivating/prototype conceptual examples before introducing the notions like \textit{$ \infty $-categories, derived schemes, higher stacks and derived stacks}.
\subsection{Derived moduli spaces} Given a site $\mathcal{C}$, which can be chosen properly according to the moduli problem of interest,  one can  define the moduli functor $\mathcal{EL}$ corresponding to a given classical field theory as 
\begin{equation}
	\mathcal{EL}: \ \mathcal{C}^{op} \longrightarrow \ Sets, \ U \mapsto \mathcal{EL}(U),
\end{equation}
where
$\mathcal{EL}(U)$ is the \textit{set} of isomorphism classes of solutions to the E-L equations over $U$.  More precisely, $\mathcal{EL}(U)$ is the moduli space $EL(U)/\mathcal{G}$ of solutions to the E-L equations modulo gauge transformation $\mathcal{G}$. But, as we discuss above, the quotient space might be problematic in general, and it fails to live in the same category. In other words, the moduli functor $\mathcal{EL}$ in general is not \emph{representable} in $\mathcal{C}.$ In order to circumvent this problem, we introduce the ``stacky" version of $\mathcal{EL}$ as the \textit{quotient stack}
\begin{equation}
	\mathcal{[EL/G]}: \ \mathcal{C}^{op} \longrightarrow  Grpds, \ \ U \mapsto \mathcal{[EL/G]}(U),
\end{equation}
where
$\mathcal{[EL/G]}(U)$ is the \textit{groupoid} of solutions to the E-L equations over $U$. Even if this explains the emergence of the stacky language in Definition \ref{main defn} in a rather intuitive way, the discussion above is just the tip of the iceberg and is still too na\"{\i}ve to capture the notion of \textit{a derived stack}.

 Thus, we need further notions in order to enjoy the richness of Definition \ref{main defn}, such as \textit{a formal neighborhood of a point} in a derived scheme/stack, a formal moduli problem, $\mathcal{L}_{\infty}$-algebras, etc... For an expository introduction to derive stacks, see \cite{Vezz2}.  

Roughly speaking, derived stacks and formal moduli problems are also higher spaces similar to stacks, but they are indeed more sensitive algebro-geometric objects to encode \textit{higher symmetries} of the theory. More precisely, let $dSt_{\mathbb{K}}$ be the $\infty$-category of derived stacks, then objects in $dSt_{\mathbb{K}}$ are simplicial presheaves preserving weak equivalences and possessing the descent/local-to-global property w.r.t. the site structure. For more details, we refer to \cite{Toen-Vezzosi}. On the other hand, \textit{formal moduli problems} are in fact particular derived stacks. 

Long story short, thanks to Yoneda's embedding, one can  realize algebro-geometric objects (like schemes, stacks, derived ``spaces", etc...) as \textit{certain functors} in addition to the standard ringed-space formulation. We have the following enlightening diagram from \cite{Vezz2} encoding such a functorial interpretation: 
\begin{center}
	\begin{tikzpicture}
		\matrix (m) [matrix of math nodes,row sep=2em,column sep=4em,minimum width=5 em] {
			CAlg_{\mathbb{K}}   & Sets  \\
			&  Grpds \\
			cdga_{\mathbb{K}}^{\leq0} &  Ssets \\};
		\path[-stealth]
		(m-1-1) edge  node [left] {{\small $\alpha$}} (m-3-1)
		edge  node [above] {{\small schemes}} (m-1-2)
		(m-1-1) edge  node [below] {} node [below] {{\small stacks}} (m-2-2)
		(m-1-1) edge  node [below] {} node [below] {{\small $ n $-stacks}} (m-3-2)
		(m-3-1) edge  node  [below] {{\small derived stacks}} (m-3-2)
		
		(m-1-2) edge  node [right] {{\small }} (m-2-2)
		(m-2-2) edge  node [right] {{\small }} (m-3-2);
	\end{tikzpicture}
\end{center} One way of interpreting this diagram is as follows: In the case of schemes, for instance, such a functorial description implies that the points of a scheme form a \textit{set}. Likewise, it implies that the collection of points of a stack has the structure of a \textit{groupoid} and not that of a set. These kinds of interpretations, in fact, suggest the name ``functor of points$ " $. So, the right hand side of the diagram in fact encodes the structure of points.

Furthermore,  since any commutative $\mathbb{K}$-algebra $A$ can be realized as an object in $cdga_\mathbb{K}^{\leq 0}$ concentrated in degree 0 with the trivial differential, the morphism $\alpha$ on the LHS of the diagram is indeed an embedding and encodes the change in the \emph{local algebraic models} of higher spaces. 

With the same spirit, the RHS of the diagram captures the level of symmetries and leads to the different ways of organizing the moduli data. I.e., the RHS is also about \emph{how to test two objects being the same}.

\paragraph{Revisiting Definition \ref{main defn}.}   Given the data of $(M, F_M, S, G)$,  where $F_M$ denotes the space of fields on some base manifold $M$, $\mathcal{S}$ is a smooth action functional on $F_M$, and $G$ is a certain symmetry group acting on $F_M$, we define a\textit{ (perturbative) classical field theory on M} to be the sheaf $\mathcal{EL}$ of derived stacks \textit{(resp. formal moduli problems)} of solutions to the E-L equations on $ M $ as follows: To each open subset $U$ of $M$, one assigns 
\begin{equation}
	U \mapsto \mathcal{EL}(U) \in dStk,
\end{equation} 
where $dStk$ denotes the $\infty$-category of derived stacks, and the derived stack $\mathcal{EL}(U)$ of solutions to the E-L equations on $ U $ is given as a functor 
\begin{equation}
	\mathcal{EL}(U): \ cdga_{k}^{\leq 0} \longrightarrow \ sSets,
\end{equation}
where $cdga_{k}^{\leq 0}$  denotes the category of commutative differential graded $k$-algebras in non-positive degrees, and $sSets$ is the $\infty$-category of simplical sets. Here $ \mathcal{EL}(U)(R) $ is the \textit{simplical set} of solutions to the defining relations (i.e. E-L equations) with values in $ R $. In other words, the points of $\mathcal{EL}(U)$ form a simplicial set.
\vspace{5pt}

As discussed above, in order to circumvent certain problems, we work with the derived moduli spaces of solutions instead of the na\"{\i}ve ones. Furthermore, we also intend to capture the \textit{perturbative} behavior of the theory, and hence this derived moduli space is defined as a \textit{formal moduli problem} \cite{Lurie}  
\begin{equation}
	\mathcal{EL}(U): \ dgArt_{k} \longrightarrow \ sSets, 
\end{equation}
where $dgArt_{k}$ denotes the category of dg artinian algebras with morphisms being simply maps of dg
commutative algebras (cf. \cite{Cosv2} Appendix A).  

\begin{remark}
	As stressed above, \textit{formal moduli problems} are in fact particular derived stacks. However, they possess a more sensitive local algebraic model comparing with the one used in derived stacks. This structure allows them to encode the \textit{perturbative natur}e of the problem under consideration in the following sense: In order to remember the pertubative behavior around the solution $p\in \mathcal{EL}(U)$, we employ the notion of a \textit{formal neighborhood} of a point (cf. \cite{Cosv2} Appendix A). This concept essentially helps us to make the underlying (higher) structure sensitive	enough to encode \textit{small thickenings} of a point obtained by adding infinitesimal directions.  To keep track of \textit{infinitesimal directions} assigned to a point $p$ (small thickenings of a point), it is indeed more suited to make use of \textit{dg artinian algebras} as  local algebraic models instead of the usual commutative $k$-algebras. That is, the underlying higher structure, informally speaking, is locally modeled on \textit{nilpotent commutative dg-algebras} such that the corresponding higher/derived space consists of points with \textit{infinitesimal directions} attached to them. 
\end{remark}
The content of this last remark may look complicated, but fortunately, formal moduli problems are tractable objects in the following way:\emph{ Every formal moduli problem can be manifested by using the language of $\mathcal{L}_{\infty}$-algebras in the sense of \cite{Lurie}}. The main ingredients of this argument will be discussed next.

\subsection{Construction of formal moduli problems} In this paper, we are interested in a particular construction, which essentially involves \textit{specific} $\mathcal{L}_{\infty}$-algebras (namely dg Lie algebras) and the corresponding Maurer-Cartan functors.

In brief, an \emph{$\mathcal{L}_{\infty}$-algebra} $\mathcal{E}$ can be considered as a generalized version of a dg Lie algebra endowed with a sequence $\{\ell_n\}$ of multilinear maps of (cohomological) degree $ 2-n $, 
\begin{equation}
	\ell_n:\mathcal{E}^{\otimes n} \longrightarrow \mathcal{E},
\end{equation}
which are called the \textit{$n$-bracket} with $n=1,2,...$ such that each bracket satisfies certain graded anti-symmetry conditions and the $n$-Jacobi rule. 

For a complete definition of an $\mathcal{L}_{\infty}$-algebra, see \cite{Cosv2}, App. A.  Notice  that  given an $\mathcal{L}_{\infty}$-algebra $(\mathcal{E}, \{\ell_n \})$, from the definition, one can observe that $\ell_1$ is indeed the differential, say $\ell_1:=d.$ 
\begin{example}
	Let $M$ be a smooth manifold and $\mathfrak{g}$ a Lie algebra. Then there exists a natural $\mathcal{L}_{\infty}$-algebra 
\begin{equation}
	\mathcal{E}:= \Omega^* (M) \otimes \mathfrak{g},
\end{equation} which will be central and appear in the context of gauge theory.
Here, the only \textit{non-zero} multilinear maps are $\ell_1:= \mathrm{d}_{dR}$ and $\ell_2:=[\cdot, \cdot ]$, where \begin{equation} \label{lie algebra structure}
	[\alpha\otimes X, \beta \otimes Y]:=\alpha\wedge\beta\otimes [X,Y]_{\mathfrak{g}}.
\end{equation} This is in fact an honest dg Lie algebra.
\end{example}

\begin{definition}
	For an $\mathcal{L}_{\infty}$-algebra $\mathcal{E} $, the Maurer-Cartan (MC) equation is 
	\begin{equation}
		\sum\limits_{n=1}^{\infty} \frac{1}{n!}\ \ell_n(\alpha^{\otimes n})=0,
	\end{equation} where $\alpha$ is an element of degree 1.
\end{definition}

\noindent Note that for the case  $ \mathcal{E}:= \Omega^* (M) \otimes \mathfrak{g} $, the MC equation becomes \begin{equation}
	\mathrm{d}A+ \frac{1}{2}[A,A]=0  \text{ for }  A\in\Omega^1 (M) \otimes \mathfrak{g}.
\end{equation}

\paragraph{Maurer-Cartan functors.}  We now wish to present the construction of a particular  formal moduli problem, which can be defined as the \textit{simplicial set of solutions to the Maurer-Cartan equation}.
\vspace{5pt}

Let $\Delta$ denote the \textit{category of finite ordered sets} with objects being finite ordered sets 
\begin{equation}
	[n]:=\{0<1<2<\cdot \cdot \cdot <n\},
\end{equation} together with the morphisms $f:[n]\rightarrow[m]$ being non-decreasing functions. Note that the set $ [n] $ corresponds to $ \Delta^n $, the usual \textit{$n$-simplex} in $\mathbb{R}^{n+1}$, given as a set\begin{equation}
	\Delta^n:= \big\{ (x_0, ... , x_n) \in \mathbb{R}^{n+1} : \sum_{i=0}^{n} x_i =1 \ \ and \ \ 0\leq x_k \leq 1 \ \ for \  all \ k \big\},
\end{equation} and hence the map $f:[n]\rightarrow[m]$ induces a linear map \begin{equation}
	f_*:\Delta^n\rightarrow\Delta^m, \ \ e_k \mapsto e_{f(k)},
\end{equation} where $e_0=(1,0,0,\cdots,0)$ and $e_k=(0,\cdots, 0, 1,0,\cdots,0)$ is the $k$th basis vector with 1 being at the $(k+1)$th slot. The ordering of $[n]$ defines a path along the edges of $ \Delta^n $,  from $e_0$ to $e_n$. For more details, see Section 1 of Appendix A in \cite{Cosv1}.

\begin{definition} \label{defn of a simplicial set}
	Let $\mathcal{C}$ be a category. A \emph{simplicial object} in $\mathcal{C}$ is a contravariant functor 
	\begin{equation}
		X_{\bullet} : \Delta \longrightarrow \mathcal{C}.
	\end{equation} 
\end{definition}If $\mathcal{C}=Sets$, then $X_{\bullet} \in Fun(\Delta,Sets)$ is called a \emph{simplicial set,} and the image $X_{\bullet} ([n])$ of $[n]$ is called the \emph{set of $ n $-simplicies} and is denoted by $ X_n. $ That is, we have \begin{equation}
	[n] \xrightarrow{X_{\bullet}}X_{\bullet}([n])=:X_n.
\end{equation} 

\begin{definition} \label{defn of MC-simplicial set}
	Let $(\mathcal{E}, \{\ell_n \})$ be an $\mathcal{L}_{\infty}$-algebra, $ (A,\mathfrak{m}) $ a dg Artinian algebra. We define the simplicial set 	$ MC(\mathcal{E}\otimes\mathfrak{m}) $ 
	of solutions to the Maurer-Cartan equation in $\mathcal{E}\otimes\mathfrak{m}$ as 
	\begin{equation}
		MC(\mathcal{E}\otimes\mathfrak{m}) \in Fun(\Delta,Sets),
	\end{equation} where an $n$-simplex in the set $MC(\mathcal{E}\otimes\mathfrak{m})_n$ of $n$-simplices is an element 
	\begin{equation} \label{element alpha}
		\alpha \in \mathcal{E}\otimes\mathfrak{m}\otimes\Omega^*(\Delta^n)
	\end{equation}  of cohomological degree 1 that satisfies the Maurer-Cartan equation.
\end{definition}

\begin{remark} \label{rmk on tensor product complex}
	In Definition \ref{defn of MC-simplicial set}, $\alpha$ is in fact an element of the tensor product complex $ \mathcal{E}\otimes\mathfrak{m}\otimes\Omega^*(\Delta^n) $ of dg modules which is defined as
	\begin{equation}
		\mathcal{E}\otimes\mathfrak{m}\otimes\Omega^*(\Delta^n)= \bigoplus \limits_k (\mathcal{E}\otimes\mathfrak{m}\otimes\Omega^*(\Delta^n))^k,
	\end{equation} 
	where $\mathcal{E}=\bigoplus_i\mathcal{E}^i$ with the differential $\mathrm{d}_{\mathcal{E}}:=\ell_1$, $\mathfrak{m}=\bigoplus _i\mathfrak{m}^i$ with the differential $\mathrm{d}_{A}$ and $\Omega^*(\Delta^n)$ is the usual de Rham complex on the $ n $-simplex $ \Delta^n $ with the de Rham differential $\mathrm{d}_{dR}$. 
	Therefore, the degree $k$ component of $ \mathcal{E}\otimes\mathfrak{m}\otimes\Omega^*(\Delta^n) $ is given by
	\begin{equation}
		(\mathcal{E}\otimes\mathfrak{m}\otimes\Omega^*(\Delta^n))^k=\bigoplus \limits_{p+q+r=k}\mathcal{E}^p\otimes\mathfrak{m}^q\otimes\Omega^r(\Delta^n),
	\end{equation}
	and hence we obtain the \textit{total complex  associated to the triple complex} 
	\begin{equation}
		\mathcal{E}\otimes\mathfrak{m}\otimes\Omega^*(\Delta^n)= \bigoplus \limits_k \bigoplus \limits_{p+q+r=k}\mathcal{E}^p\otimes\mathfrak{m}^q\otimes\Omega^r(\Delta^n),
	\end{equation} with the \textit{total differential} $ d_{tot}^k: (\mathcal{E}\otimes\mathfrak{m}\otimes\Omega^*(\Delta^n))^k \rightarrow (\mathcal{E}\otimes\mathfrak{m}\otimes\Omega^*(\Delta^n))^{k+1}$ defined by \begin{equation}
		d_{tot}^k=\sum \limits_{p+q+r=k} d_1^{p,q,r}+(-1)^pd_2^{p,q,r}+(-1)^{p+q}d_3^{p,q,r}, 
	\end{equation} where
	\begin{align}
		d_1^{p,q,r} &= \mathrm{d}_{\mathcal{E}}^p \otimes id_{A}^q\otimes id_{\Omega^r}, \ \ \mathrm{d}_{\mathcal{E}}^p: \mathcal{E}^p \rightarrow \mathcal{E}^{p+1} \nonumber \\
		d_2^{p,q,r} &= id_{\mathcal{E}}^p \otimes \mathrm{d}_{A}^q\otimes id_{\Omega^r}, \ \ \mathrm{d}_{A}^q: \mathfrak{m}^q \rightarrow \mathfrak{m}^{q+1} \nonumber \\
		d_3^{p,q,r} &= id_{\mathcal{E}}^p \otimes id_{A}^q\otimes \mathrm{d}_{dR}^r, \ \ \mathrm{d}_{dR}^r: \Omega^r \rightarrow \Omega^{r+1}. 
	\end{align}
	For a more concrete treatment to the notions like double/triple complexes and their total complexes, see Chapter 12 of \cite{Stacks}.  In order to illustrate the situation related to the triple complexes and motivate the structure of such ``higher dimensional" cochain complexes, we examine the following simplest case. 
	
	Suppose that $A$ is an ordinary $\mathrm{k}$-algebra, and $\mathfrak{m}:=\mathrm{k}$. Note that $A$ can be viewed as a complex that is concentrated at degree 0, and all other components are trivial with differential being zero. Hence, in this situation, we can consider  $ \mathcal{E}\otimes\mathfrak{m}\otimes\Omega^*(\Delta^n) $ as a double complex and write $ \mathcal{E}\otimes\Omega^*(\Delta^n)$ instead. Furthermore, we diagrammatically have

	\begin{equation} \label{double cmpx diagram}
		\begin{tikzpicture}
			\matrix (m) [matrix of math nodes,row sep=2em,column sep=3em,minimum width=2em] 
			{ 
				& \vdots         & \vdots &      \\
				\cdots  & \mathcal{E}^p\otimes\Omega^{r+1}(\Delta^n) &\mathcal{E}^{p+1}\otimes\Omega^{r+1}(\Delta^n)       &\cdots\\
				\cdots  & \mathcal{E}^p\otimes\Omega^{r}(\Delta^n) & \mathcal{E}^{p+1}\otimes\Omega^{r}(\Delta^n)      &\cdots \\
				& \vdots         & \vdots &        \\
			};
			\path[-stealth]
			(m-2-2) edge  node [left] {{\small $  $}} (m-1-2)
			(m-3-2) edge  node [left] {{\small $ d_2^{p,r} $}} (m-2-2)
			(m-4-2) edge  node [left] {{\small $  $}} (m-3-2)
			(m-2-1.east|-m-2-2) edge  node [below] {} node [below] {{\small $  $}} (m-2-2)
			(m-3-1.east|-m-3-2) edge  node [below] {} node [below] {{\small $  $}} (m-3-2)
			
			(m-2-2.east|-m-2-3) edge  node [below] {} node [above] {{\small $ d_1^{p,r+1} $}} (m-2-3)
			(m-3-2.east|-m-3-3) edge  node [below] {} node [below] {{\small $ d_1^{p,r} $}} (m-3-3)
			
			(m-2-3.east|-m-2-4) edge  node [below] {} node [below] {{\small $  $}} (m-2-4)
			(m-3-3.east|-m-3-4) edge  node [below] {} node [below] {{\small $  $}} (m-3-4)
			
			(m-2-3) edge  node [right] {{\small $ $}} (m-1-3)
			(m-3-3) edge  node [right] {{\small $ d_2^{p+1,r}$}} (m-2-3)
			(m-4-3) edge  node [right] {{\small $ $}} (m-3-3);
			
		\end{tikzpicture}
	\end{equation}
	where $ d_1^{p,r} = \mathrm{d}_{\mathcal{E}}^p \otimes id_{\Omega^r} $ and $ d_2^{p,r} = id_{\mathcal{E}}^p \otimes \mathrm{d}_{dR}^r $ for all $p,r$. Note that each square in the diagram is commutative, and hence different parts of the differential are compatible. For the precise structural relations, we again refer to Ch. 12 of \cite{Stacks}.
\end{remark}

\begin{definition} \label{defn of MC functor Bg}
	Given an $\mathcal{L}_{\infty}$-algebra $\mathcal{E}$, we  define the \emph{Maurer-Cartan functor} $\mathcal{MC}_{\mathcal{E}}$ for  $\mathcal{E}$ by 
	\begin{equation}
		\mathcal{MC}_{\mathcal{E}}: \ dgArt_{k} \longrightarrow \ sSets, \ \ (A,\mathfrak{m})\longmapsto \mathcal{MC}_{\mathcal{E}}[(A,\mathfrak{m})]:= MC(\mathcal{E}\otimes\mathfrak{m}), 
	\end{equation} where the set of $n$-simplicies of $MC(\mathcal{E}\otimes\mathfrak{m})$ is defined as above (Definition \ref{defn of MC-simplicial set}):
	\begin{equation} \label{defn of the set of n-simplicies of MC-simplicial set}
		MC(\mathcal{E}\otimes\mathfrak{m})_n = \Bigg\{	\alpha \in \bigoplus \limits_{p+q+r=1}\mathcal{E}^p\otimes\mathfrak{m}^q\otimes\Omega^r(\Delta^n) :	\	\mathrm{d}\alpha+\sum\limits_{n=2}^{\infty} \frac{1}{n!}\ell_n(\alpha^{\otimes n})=0		\Bigg\}
	\end{equation}
\end{definition}

\begin{lemma}
	The  functor $\mathcal{MC}_{\mathcal{E}}$ is a formal moduli problem. 
\end{lemma} \noindent For more details on the construction of the  functor $\mathcal{MC}_{\mathcal{E}}$, we refer to Section 4.1 of \cite{Cosv2}. 
\begin{theorem} (\cite{Lurie}, Theorem 2.0.2) \label{Lurie s thm}
	Every formal moduli problem is represented by a \textit{Maurer-Cartan functor} $\mathcal{MC}_{\mathcal{E}}$ for some differential graded Lie algebra (or an $\mathcal{L}_{\infty}$ algebra) $\mathcal{E} $ up to a weak equivalence. More precisely, there exists an equivalence of $\infty$-categories \begin{equation}
		dgla_k\xrightarrow{\sim}Moduli_k \subset Fun(dgArt_k,Ssets),
	\end{equation} where $ dgla_k $ and $Moduli_k$ denote $\infty$-categories of differential graded Lie algebras over $k$ and that of formal moduli problems over $k$ respectively with $k$ being a field of characteristic zero. 
\end{theorem}

\begin{remark}
	Given a pair of dg Lie algebras $(g,d)$ and $(g',d')$, \textit{a map of dg Lie algebras} from $(g,d)$ to $(g',d')$ is a map of cochain complexes $f: (g,d) \rightarrow (g',d')$ s.t. $f([X,Y]_g)=[f(X),f(Y)]_{g'}$ for $X \in g_i$ and $Y\in g_j$. Also, we will say that a map of dg Lie algebras $f: (g,d) \rightarrow (g',d')$ is a \textit{quasi-isomorphism} if the underlying map of cochain complexes is a quasi-isomorphism. Here, it follows from Proposition 2.1.10 in \cite{Lurie} that $ dgla_k $ admits the structure of a model category with weak equivalences being quasi-isomorphisms. In that respect, two dg Lie algebras $\mathcal{E}$ and $\mathcal{E}'$ induce equivalent formal moduli problems only if they can be joined by a chain of quasi-isomorphisms. In particular, every dg Lie algebra over an algebraically  closed field $k$ of characteristic zero determines a formal moduli problem. For more details, we refer to \cite{Lurie}.
\end{remark}

\section {Final remarks and recasting some examples} \label{sections_recasting some examples}
Now we are in place of summarizing what we have done so far and provide some key observations leading to the derived reformulation of classical field theories: 
\begin{enumerate}
	\item[\textit{i.}] Describing a classical field theory is equivalent to a study of the moduli space $\mathcal{EL}$ of solutions to the Euler-Lagrange equations (and hence the critical locus of action functional), which in fact corresponds to a certain moduli functor.
	\vspace{5pt}
	\item[\textit{ii.}]  As stressed before, a typical moduli functor, however, is not representable in general due to certain problems, such as the existence of degenerate critical points or non-freeness of the action of the symmetry group acting on the space of fields. In order to avoid problems of these kinds (and to capture the perturbative behavior at the same time), one requires to adopt the language of derived algebraic geometry, and hence one needs to replace the na\"{\i}ve notion of moduli problem by the so-called \textit{formal moduli problem} in the sense of Lurie as discussed above.  
	\vspace{5pt}
	\item[\textit{iii.}] Formal moduli problems $\mathcal{F}$, on the other hand, are unexpectedly tractable notions (thanks to Lurie's result, Theorem \ref{Lurie s thm}) in the sense that understanding $\mathcal{F}$, at the end of the day, boils down to finding a suitable (local) dgla (or $\mathcal{L}_{\infty}$-algebra) $\mathcal{E}$ such that $\mathcal{F}$ can be represented by the Maurer-Cartan functor $ \mathcal{MC}_{\mathcal{E} } $ associated to $\mathcal{E} $.
	\vspace{5pt}
	\item[\textit{iv.}] Having obtained appropriate dgla (or $\mathcal{L}_{\infty}$-algebra) $\mathcal{E}$, one needs to analyze the structure of $\mathcal{E}$ so as to encode the theory under consideration.
	
\end{enumerate}

\noindent\textbf{Revisiting Example \ref{example_free scalar massive field theory}.} Consider a free scalar \textit{massless} field theory on a Riemannian manifold $M$ with space of field being $ C^{\infty}(M)$  and the action functional governing the theory as 
\begin{equation}
	\mathcal{S}(\phi):= \int\limits_{M} \phi \Delta \phi.
\end{equation} 
The corresponding E-L equation in this case turns out to be 
\begin{equation}
	\Delta\phi=0,
\end{equation} and hence the moduli space $EL$ of solutions to the E-L equations is the moduli space of harmonic functions \begin{equation}
	\big\{ \phi \in  C^{\infty}(M) :   \Delta\phi=0        \big\}.
\end{equation} Now, having employed the derived enrichment $\mathcal{EL}$  of $ EL $ as described above, we need to find a suitable $\mathcal{L}_{\infty}$-algebra $\mathcal{E} $ whose Maurer-Cartan fuctor $ \mathcal{MC}_{\mathcal{E} } $ represents the formal moduli problem $\mathcal{EL}$. The answer is as follows: We define $ \mathcal{E} $ to be the two-term complex 
\begin{equation}
	\mathcal{E}: C^{\infty}(M)\xrightarrow{\Delta}C^{\infty}(M)[-1],\end{equation} which is concentrated in degree 0 and 1 with a sequence $\{\ell_n\}$ of multilinear maps, where $\ell_1:=\Delta$ and $\ell_i=0$ for all $i\neq1.$  The Maurer-Cartan equation, on the other hand, turns out to be \begin{equation}
	\Delta\phi=0,
\end{equation}and hence the set  of 0-simplices of the simplicial set $ \mathcal{MC}_{\mathcal{E} }(A) $ for $A$ ordinary Artinian algebra is given as \begin{equation}
	\big\{ \phi :   \Delta\phi=0        \big\},
\end{equation} as desired.
For further details and interpretation of other simplices, see Chapter 2 of \cite{Cosv1} or Chapter 4 of \cite{Cosv2}.

\vspace{5pt}

\noindent\textbf{Revisiting Example \ref{example_CS theory}.}  We shall revisit $SU(2)$ Chern-Simons gauge theory on a closed, orientable 3-manifold $ X $. As before, Let $P\rightarrow X$ be a trivial principal $SU(2)$-bundle on $X$, its Lie algebra $ \mathfrak{g}:=\mathfrak{su}(2)$, and $A\in \mathcal{A}:=\Omega^1 (X, \mathfrak{su}(2))$ the Lie algebra-valued connection 1-form on $X$ together with the Chern-Simons action funtional $ CS: \mathcal{A} \longrightarrow S^1$ given by
\begin{equation}
	CS(A):=\frac{k}{4\pi} \displaystyle \int \limits_{X} \mathrm{Tr}(A\wedge \mathrm{d} A +\frac{2}{3}  A \wedge A \wedge A), ~~~~ k\in \mathbb{Z}.
\end{equation}Here the gauge group $\mathcal{G}=Map(X,SU(2))$ acting on the space $\mathcal{A}$ as usual.
The corresponding E-L equation in this case turns out to be 
\begin{equation}
	F_{A}:=\mathrm{d} A+A \wedge A=0,
\end{equation}
where $F_{A}$ is the curvature two-form on $X$ associated to $A$. 

Now, we are interested in the moduli space $ \mathcal{M}_{flat}$ of flat connections modulo gauge transformations \begin{equation}
	\{ A\in  \Omega^1 (X, \mathfrak{su}(2)) :  F_A=0        \} \big/ \sim.
\end{equation}
As before, having introduced  derived counterpart $\mathcal{EL}_{flat}$ of $ \mathcal{M}_{flat}$, we define a suitable $\mathcal{L}_{\infty}$-algebra $\mathcal{E} $ encoding the formal moduli problem as follows:
\begin{equation}
	\mathcal{E}:= \Omega^* (X) \otimes \mathfrak{g}[1],
\end{equation}
where the only \textit{non-zero} multilinear maps are $\ell_1:= \mathrm{d}_{dR}$ and $\ell_2:=[\cdot, \cdot ]$ given as in Equation \ref{lie algebra structure}. Notice that the Maurer-Cartan equation in this case becomes \begin{equation}
	\mathrm{d}A+ \frac{1}{2}[A,A] =0,
\end{equation} and hence the corresponding the Maurer-Cartan functor $ \mathcal{MC}_{\mathcal{E} } $ yields the desired result. For a complete treatment, you may check Chapter 4 of \cite{Cosv2}, Chapter 5.4 of \cite{Costello Renormalizetion book}, or \cite{Costello Renormalisation and BV formalism paper}. 

Furthermore, as stressed in Chapter 5.4 of \cite{Costello Renormalizetion book}, the spaces of all fields associated to the theory are encoded by the $\mathcal{L}_{\infty}$-algebra $\mathcal{E} $ in the Batalin-Vilkovisky formalism as follows:

\begin{itemize}
	\item The space of degree $ -1 $ fields, called \textit{ghosts}, corresponds to the space \begin{equation}
		\Omega^0 (X) \otimes \mathfrak{g}=Map(X,SU(2))=\mathcal{G}.
	\end{equation}
	\item The space of degree $ 0 $ fields, called \textit{fields}, corresponds to the space  \begin{equation}
		\Omega^1 (X) \otimes \mathfrak{g}.
	\end{equation} 
	\item The space of degree $ 1 $ fields, called \textit{anti-fields}, corresponds to the space  \begin{equation}
		\Omega^2 (X) \otimes \mathfrak{g}.
	\end{equation} 
	\item The space of degree $ 2 $ fields, called \textit{anti-ghosts}, corresponds to the space  \begin{equation}
		\Omega^3 (X) \otimes \mathfrak{g}.
	\end{equation} 
\end{itemize}

\section*{Acknowledgments}

This note serves as an  survey on  derived-geometric structures in physics. The main body of the text contains key ideas and results that are gathered from the literature. Arguments we present here are very well-known to the experts, and as a disclaimer they are not meant to provide neither original nor new results (nor a complete reference list) related to either of the subjects mentioned above. But we hope that the material we present in this paper can provide a brief introduction and a na\"{\i}ve guideline to the existing literature for non-experts who may wish to learn the subject. So, one may propose that our systematic presentations of the various subjects may come in handy. 

I am very grateful to Ali Ulaş Özgür Kişisel and Bayram Tekin for their enlightening, fruitful and enjoyable conversations during our regular research meetings that essentially lead to the preparation of this note.  

I would also like to thank the associate editor and anonymous reviewers for their valuable comments and suggestions which improved and clarified the manuscript.




\begin{thebibliography}{0}
\bibitem{Mnev} P. Mnev,  Lectures on BV formalism and its applications in topological quantum field theory, math-ph/1707.08096v1.
	
	\newpage
\end{thebibliography}
\end{document}